%% file: paper_v3.tex
\documentclass[]{article}
\usepackage{proceed2e}

\usepackage{times}

\newcommand*{\rd}{\mathrm{d}}

\input{mypreamble.tex}
\title{Active Uncertainty Calibration in Bayesian ODE Solvers}

\author{ {\bf Hans Kersting and Philipp Hennig} \\
Max-Planck-Institute for Intelligent Systems \\
Spemannstra{\ss}e, 72076 T{\"u}bingen, Germany \\
\texttt{[hkersting|phennig]@tue.mpg.de}
}

\begin{document}

\maketitle

\begin{abstract}
There is resurging interest, in statistics and machine learning, in solvers for ordinary differential equations (ODEs) that return probability measures instead of point estimates. Recently, Conrad et al.~introduced a sampling-based class of methods that are `well-calibrated' in a specific sense. But the computational cost of these methods is significantly above that of classic methods. On the other hand, Schober et al.~pointed out a precise connection between classic Runge--Kutta ODE solvers and Gaussian filters, which gives only a rough probabilistic calibration, but at negligible cost overhead. By formulating the solution of ODEs as approximate inference in linear Gaussian SDEs, we investigate a range of probabilistic ODE solvers, that bridge the trade-off between computational cost and probabilistic calibration, and identify the inaccurate gradient measurement as a crucial source of uncertainty. We propose the novel filtering-based method Bayesian Quadrature filtering (BQF) which uses Bayesian quadrature to actively learn the imprecision in the gradient measurement by collecting multiple gradient evaluations.
\end{abstract}

\section{INTRODUCTION}

The numerical solution of an initial value problem (IVP) based on an \emph{ordinary differential equation} (ODE)
\begin{equation} \label{orODE}
  u^{(n)}(t) = f(t,u(t),\dots,u^{(n-1)}(t)),\ u(0)=u_0 \in \R^D,
\end{equation}
of order $n\in\N$, with $u:\mathbb R\to \mathbb R^D$, $f: [0,T] \times \mathbb R^{nD}\to \mathbb R^D$, $T>0$, is an essential topic of numerical mathematics, because ODEs are the standard model for dynamical systems. Solving ODEs with initial values is an exceedingly well-studied problem \citep[see][for a comprehensive presentation]{hairer87:_solvin_ordin_differ_equat_i} and modern solvers are designed very efficiently. Usually, the original ODE \eqref{orODE} of order $n$ is reduced to a system of $n$ ODEs of first order
\begin{equation}
  u^{\prime}(t) = f(t,u(t)),\; u(0)=u_0 \in \R^D,
\end{equation}
which are solved individually. The most popular solvers in practice are based on some form of Runge--Kutta (RK) method (as first introduced in \citet{Runge} and \citet{Kutta}) which employ a  weighted sum of a fixed amount of gradients in order to iteratively extrapolate a discretized solution. That is, these methods collect `observations' of approximate gradients of the solved ODE, by evaluating the vector field $f$ at an estimated solution, which is a linear combination of previously collected `observations':
\begin{equation}
  y_i = f\left(t+c_ih,u_0+\sum_{j<i} w_{ij} y_j\right).
\end{equation}
The final extrapolation step is a weighted sum of these gradients:
\begin{equation}
  \hat u(t+h) = u(t) + \sum_{i<s} b_i y_i.
\end{equation}
The weights of $s$-stage RK methods of $p$-th order are carefully chosen so that the numerical approximation $\hat{u}$ and the Taylor series of the exact solution $u$ coincide up to the term $h^p$, thereby yielding a local truncation error of high polynomial order,
\begin{equation}\| u(t_0 + h) - \hat u(t_0+h)\|=\mathcal O(h^{p+1}),\end{equation} 
for $h\to 0$. One can prove that $s\geq p$ in general, but for $p\leq 4$ there are RK methods with $p=s$. Hence, allowing for more function evaluations can drastically improve the speed of convergence to the exact solution. 

The polynomial convergence is impressive and helpful; but it does not actually quantify the inevitable epistemic uncertainty over the accuracy of the approximate solution $\hat{u}$ for a concrete non-vanishing step-size $h$. One reason one may be concerned about this in machine learning is that ODEs are often one link of a chain of algorithms performing some statistical analysis. When employing classic ODE solvers and just plugging in the solution of the numerical methods in subsequent steps, the resulting uncertainty of the whole computation is ill-founded, resulting in overconfidence in a possibly wrong solution. It is thus desirable to model the epistemic uncertainty. Probability theory provides the framework to do so. Meaningful probability measures of the uncertainty about the result of deterministic computations (such as ODE solvers) can then be combined with probability measures modeling other sources of uncertainty, including `real' aleatoric randomness (from e.g.~sampling). Apart from quantifying our certainty over a computation, pinning down the main sources of uncertainty could furthermore improve the numerical solution and facilitate a more efficient allocation of the limited computational budget.

A closed framework to measure uncertainty over numerical computations was proposed by \citet{skilling1991bayesian} who pointed out that numerical methods can be recast as statistical inference of the latent exact solution based on the observable results of tractable computations. 
In this spirit, \citet{HennigAISTATS2014} phrased this notion more formally, as Gaussian process (GP) regression. Their algorithm class, however, could not guarantee the high polynomial convergence orders of Runge--Kutta methods. 
In parallel development, \citet{chkrebtii_old_arxiv} also introduced a probabilistic ODE solver of similar structure (i.e.~based on a GP model), but using a Monte Carlo updating scheme. These authors showed a linear convergence rate of their solver, but again not the high-order convergence of classic solvers. \\
Recently, \citet{schober2014nips} solved this problem by finding prior covariance functions which produce GP ODE solvers whose posterior means \emph{exactly} match those of the optimal Runge--Kutta families of first, second and third order.
While producing only a slight computational overhead compared to classic Runge--Kutta, this algorithm---as any GP-based algorithm---only returns Gaussian measures over the solution space. \\
In contrast, \citet{conrad_old_arxiv} recently provided a novel sampling-based class of ODE solvers which returns flexible non-Gaussian measures over the solution space, but creates significant computational overhead by running the whole classic ODE solvers multiple times over the whole time interval $[0,T]$ in order to obtain meaningful approximations for the desired measure. \\
For practitioners, there is a trade-off between the desire for quantified uncertainty on the one hand, and low computational cost on the other. The currently available probabilistic solvers for ODEs either provide only a roughly calibrated uncertainty \citep{schober2014nips} at negligible overhead or a more fine-grained uncertainty supported by theoretical analysis \citep{conrad_old_arxiv}, at a computational cost increase so high that it rules out most practical applications.
In an attempt to remedy this problem, we propose an algorithm enhancing the method of \citet{schober2014nips} by improving the gradient measurement using modern probabilistic integration methods. By modeling the  uncertainty where it arises, i.e.~the imprecise prediction of where to evaluate $f$, we hope to gain better knowledge of the propagated uncertainty and arrive at well-calibrated posterior variances as uncertainty measures.

\section{BACKGROUND}\label{Background}



\subsection{SAMPLING-BASED ODE SOLVERS} 
\label{Sampling}


The probabilistic ODE solver by \citet{conrad_old_arxiv} modifies a classic deterministic one-step numerical integrator $\Psi_h$ (e.g. Runge--Kutta or multiderivative methods, cf. \citet{hairer87:_solvin_ordin_differ_equat_i}) and models the discretization error of $\Psi_h$ by adding suitably scaled i.i.d. Gaussian random variables $\{\xi_k\}_{k=0,\dots,K}$ after every step. Hence, it returns a discrete solution $\{U_k\}_{k=0,\dots,K}$ on a mesh $\{t_k = kh \}_{k=0,\dots,K}$ according to the rule
\begin{equation}
  U_{k+1} = \Psi_h (U_k) + \xi_k. \label{Conrad_algo}
\end{equation}
This discrete solution can be extended into a continuous time approximation of the ODE, which is random by construction and can therefore be interpreted as a draw from a non-parametric probability measure $Q_h$ on the solution space $C^1\left([0,T],\R^n\right)$, the Banach space of continuously differentiable functions. This probability measure can then be interpreted as a notion of epistemic uncertainty about the solution. This is correct in so far as, under suitable assumptions, including a bound on the variance of the Gaussian noise, the method converges to the exact solution, in the sense that $Q_h$ contracts to the Dirac measure on the exact solution $\del_u$ with the \emph{same} convergence rate as the original numerical integrator $\Psi_h$, for $h\to 0$:
If $(\xi_{k,h})_{k=1}^N \sim \mathcal N(0,\Var(h))$ with $\Var (h) = \mathcal{O}(h^{2q+1})$, then
  \begin{equation}
    \sup_{0\leq kh\leq T} \E^h \left \Vert u_k - U_k \right \Vert^2 \leq \sigma\cdot h^{2q}.
  \end{equation}
This is a significant step towards a well-founded notion of uncertainty calibration for ODE solvers: It provides a probabilistic extension to classic method which does not break the convergence rate of these methods.

In practice, however, the precise shape of $Q_h$ is \emph{not} known and $Q_h$ can only be interrogated by sampling, i.e.~repeatedly running the entire probabilistic solver. After $S$ samples, $Q_h$ can be approximated by an empirical measure $Q_h(S)$. In particular, the estimated solution and uncertainty can only be expressed in terms of statistics of $Q_h(S)$, e.g.~by the usual choices of the empirical mean and empirical variance respectively or alternatively by confidence intervals. For $S\to\infty$, $Q_h(S)$ converges in distribution to $Q_h$ which again converges in distribution to $\del_u$ for $h\to 0$: 
\begin{equation} Q_h(S) \stackrel{S\to\infty}{\to} Q_h \stackrel{h\to 0}{\to} \del_u. \label{convergence}
\end{equation}
The theoretical mathematics in \citet{conrad_old_arxiv} only concerns the convergence of the latent probability measures $\{Q_h\}_{h>0}$. Only the empirical measures $\{Q_h(S)\}_{S\in\N}$, however, can be observed. Consequently, it remains unclear whether the empirical mean of $Q_h(S)$ for a fixed step-size $h>0$ converges to the exact solution as $S\to\infty$ and whether the empirical variance of $Q_h(S)$ is directly related, in an analytical sense, to the approximation error. In order to extend the given convergence results to the practically observable measures $\{Q_h(S)\}_{S\in\N}$ an analysis of the first convergence in~\eqref{convergence} remains missing. The deterministic algorithm proposed below avoids this problem, by instead constructing a (locally parametric) measure from prior assumptions.

The computational cost of this method also seems to mainly depend on the rate of convergence of $Q_h(S) \to Q_h$ which determines how many (possibly expensive) runs of the numerical integrator $\Psi_h$ over $[0,T]$ have to be computed and how many samples have to be stored for a sufficient approximation of $Q_h$. 
Furthermore, we expect that in practice the mean of $Q_h$, as approximated by $Q_h(S)$ might not be the best possible approximation, since in one step the random perturbation of the predicted solution by Gaussian noise $\xi_k$ worsens our solution estimate with a probability of more than $1/2$, since---due to the full support of Gaussian distributions---the numerical sample solution is as likely to be perturbed away from as towards the exact solution and---due to the tails of Gaussian distributions---it can also be perturbed way past the exact solution with positive probability.

\subsection{A FRAMEWORK FOR GAUSSIAN FILTERING FOR ODES}\label{Framework}


Describing the solution of ODEs as inference in a joint Gaussian model leverages state-space structure to achieve efficient inference. Therefore, we employ a Gauss--Markov prior on the state-space: \emph{A priori} we model the solution function and $(q-1)$ derivatives $( u, \dot{u}, u^{(2)},\dots,u^{(q-1)} ):[0,T] \to \mathbb \R^{qD}$ as a draw from a $q$-times integrated Wiener process $X=(X_t)_{t\in [0,T]} = (X_t^{(1)},\dots,X_t^{(q)})^{T}_{t\in [0,T]}$, i.e.~the dynamics of $X_t$ are given by the linear stochastic differential equation \citep{karatzas1991brownian,oksendal2003stochastic}:
\begin{align}\label{SDE}
dX_t &= F X_t dt +  L dW_t,  \\
X_0 &= \xi,\quad \xi \sim \mathcal N(m(0),P(0)),
\end{align}
with constant drift $F\in \R^{q\times q}$ and diffusion $L \in \R^{q}$ given by
\begin{align} 
    F = \begin{pmatrix} 0 & f_1 & 0 & \hdots & 0 \\ 0 & 0 & f_2 & \hdots & 0 \\ \vdots & & \ddots & \ddots & \vdots \\ 0 & & \hdots & 0 & f_{q-1} \\ 0 & & \hdots & 0 & 0     \end{pmatrix},\; L = \begin{pmatrix} 0 \\ 0 \\ \vdots \\ 0 \\ \sigma   \end{pmatrix}
\end{align}
for all $t\in [0,T]$ and some $f_1,\dots,f_{q-1}\in\R$, where $W_t$ denotes a $q$-dimensional Wiener process ($q\geq n$). Hence, we are a priori expecting that $u^{(q)}$ behaves like a Brownian motion with variance $\sigma^2$ and that $u^{(i)}$ is modeled by $(q-1-i)$-times integrating this Brownian motion. The fact that the $(i+1)$-th component is the derivative of the $i$-th component in our state space is captured by a drift matrix with non-zero entries only on the first off-diagonal. The entries $ f_1,\dots, f_{q-1}$ are damping factors. A standard choice is e.g.~$ f_i = i$. Without additional information, it seems natural to put white noise on the $q$-th derivative as the first derivative which is not captured in the state space. This gives rise to Brownian noise on the $(q-1)$-th derivative which is encoded in the diffusion 
matrix scaled by variance $\sigma^2$. Hence, we consider the integrated Wiener process a natural prior. For notational simplicity, only the case of scalar-valued functions, i.e.~$D=1$, is presented in the following. The framework can be extended to $D\geq 2$ in a straightforward way by modeling the output dimensions of $f$ as independent stochastic processes.

Since $X$ is the strong solution of a linear equation \eqref{SDE} with normally distributed initial value $X_0$, it follows from the theory of linear SDEs \citep{karatzas1991brownian} that $X$ is a uniquely-determined Gauss--Markov process. This enables Bayesian inference in a highly efficient way by Gaussian filtering \citep{saatci2011thesis}). For time invariant linear SDEs like \eqref{SDE}, the fixed matrices for Gaussian filtering can be precomputed analytically \citep{sarkka2006thesis}.


In addition, \citet{schober2014nips} showed that for $q\leq 3$ inference in this linear SDE yields Runge--Kutta steps.

Equipped with this advantageous prior we can perform Bayesian inference. The linearity and time-invariance of the underlying SDE permits to formulate the computation of the posterior as a Kalman filter (KF) (cf. \citep{Sarkka2013} for a comprehensive introduction) with step size $h>0$. The prediction step of the KF is given by
\begin{align}
  m_{t+h}^- &= A(h)m_{t}, \label{KFP1}  \\
  P_{t+h}^- &= A(h)P_tA(h)^{T}+Q(h), \label{KFP2}
\end{align}
with matrices $A(h), Q(h) \in \R^{q \times q }$  with entries
\begin{align}
  A(h)_{i,j} = &\exp (h F)_{i,j} = \chi_{j\geq i} \frac {h^{j-i}}{(j-i)!} \left ( \prod_{k=0}^{j-i-1} f_{i+k} \right ),   
  \notag
  \\
  Q(h)_{i,j} = &\sigma^2 \left( \prod_{k_1=0}^{q-1-i} f_{i+k_1} \right) \cdot \left( \prod_{k_2=0}^{q-1-j} f_{j+k_2} \right)\cdot \notag \\ &\frac{h^{2q+1-i-j}}{(q-i)!(q-j)!(2q+1-i-j)}.
\end{align}
It is followed by the update step
\begin{align} 
  z &= y - H m_{t+h}^-, \label{KFU1}  \\
  S &= H P_{t+h}^- H^{T} + R, \label{KFU2}  \\
  K &= P_{t+h}^- H^{T} S^{-1}, \label{KFU3} \\
  m_{t+h} &= m_{t+h}^- + Kz, \label{KFU4} \\
  P_{t+h} &= P_{t+h}^- - KHP_{t+h}^-, \label{KFU5} 
\end{align}
where $ H = e_{n}^{T} \in \R^{1\times q} $ is the $n$-th unit vector.  

Between the prediction and update step the $n$-th derivative of the exact solution $\frac{\partial^{n} u}{\partial x^{n}}$ at time $t+h$ as a measurement for the $n$-th derivative and the noise of this measurement are estimated by the variable $y$ and $R$. In order to derive precise values of $y$ and $R$ from the Gaussian prediction $\mathcal N(m_{t+h},P_{t+h})$, we would have to compute the integrals
\begin{equation}
  y = \int f(t+h,m_{t+h}^- + x) \mathcal N(x;0,P^-_{t+h})\ \rd x  
\end{equation}
and
\begin{align} R = \int &f(t+h,m_{t+h}^- + x) f(t+h,m_{t+h}^- + x)^T \cdot \notag  \\ &\mathcal N (x;0,P^-_{t+h}) \ \rd x -yy^T, \end{align}
which are intractable for most choices of $f$. Below we investigate different ways to address the challenge of accurately approximating these integrals while not creating too much computational overhead.

\subsection{MEASUREMENT GENERATION OPTIONS FOR GAUSSIAN FILTERING}\label{Measurement Generation}

\citet{schober2014nips} as, to the best of our knowledge, the first ones to point out the connection between Gaussian filtering and probabilistic ODE solvers, presents an algorithm which simply evaluates the gradient at the predicted mean, which is equivalent to setting $y$ to be equal to its maximum likelihood estimator: \begin{equation} y = f(t+h,m_{t+h}^-),\; R=0. \end{equation}

While ensuring maximum speed, this is clearly not an ideal measurement. In our atomless predicted probability measure $\mathcal{N}(m_{t+h}^-,P_{t+h}^-)$ the mean predictor $m_{t+h}^-$ is different from its exact value $(u^{(0)}(t+h),\dots,u^{(n)}(t+h))^T$ almost surely. Hence, for a non-constant $f$ the estimate will be inaccurate most of the times. In particular this method deals poorly with `skewed' gradient fields (a problem that leads to a phenomenon known as `Lady Windermeres fan' \citep{hairer87:_solvin_ordin_differ_equat_i}). To get a better estimate of the exact value of $y$, more evaluations of $f$ seem necessary.

Therefore, we want to find numerical integration methods which capture $y$ and $R$ with sufficient precision, while using a minimal number of evaluations of $f$. Possible choices are:

\begin{enumerate}[(i)]

\item \emph{Monte Carlo integration by sampling}: 
\begin{align}
  y &= \frac{1}{N} \sum_{i=1}^N f(t+h,x_i), \label{MCIS1}  \\ R &= \frac 1N \sum_{i=1}^N f(t+h,x_i) f(t+h,x_i)^T - y y^T, \label{MCIS2} \\ x_i &\sim \mathcal N( m^-_{t+h}, P^-_{t+h} ),  \label{MCIS3}
\end{align}

(which is \emph{not} the same as the sampling over the whole time axis in \citep{conrad_old_arxiv}).
\item \emph{Approximation by a first-order Taylor series expansion:} \begin{align} f&(t+h,m_{t+h}^- + x)\notag \\ &\simeq f(t+h,m_{t+h}^-) + \nabla f(t+h,m_{t+h}^- + x)\cdot x \end{align} and thereby deriving moments of the linear transform of Gaussian distributions:  
\begin{align}
  y &= f(t+h,m_{t+h}^-),  \\
  R &= \nabla f(t+h,m_{t+h}^-) P_{t+h}^- \nabla f(t+h,m_{t+h}^-)^T.
\end{align}
\item Integration by \emph{Bayesian quadrature} with Gaussian weight function: 
\begin{align}
  y &= \alpha^T K^{-1} \begin{pmatrix} f(x_1) , \hdots , f(x_n) \end{pmatrix}^T,  \label{BQy}  \\
  R &= \int \int k(x,x^{\prime}) w(x) w(x^{\prime}) \ \rd x\rd x^{\prime} - \alpha^T K^{-1} \alpha.  \label{BQR}
\end{align}
with $w(x) = \mathcal N(x;m^-_{t+h},P^-_{t+h})$, kernel matrix $K \in \R^{N \times N}$ with $K_{i,j} = k(x_i,x_j)$ and $\alpha =(\alpha(1),\dots,\alpha(N))^T \in \R^N$ with $\alpha (i) = \int k(x,x_i) w(x)\ \rd x $ for a predefined covariance function $k$ and evaluation points $(x_i)_{i=1,\dots,N}$ (cf.~section \ref{BQ}).
\end{enumerate}
Our experiments, presented in Section \ref{experiments}, suggest that BQ is the most useful option.

Monte Carlo integration by sampling behaves poorly if the trajectory of the numerical solution passes through domain areas (as e.g.~in the spikes of oscillators governed by non-stiff ODEs) where $f$ takes highly volatile values since the random spread of samples from the domain are likely to return a skewed spread of values resulting in bad predictions of $y$ with huge uncertainty $R$. Hence, the posterior variance explodes and the mean drifts back to its zero prior mean, i.e. $m_{t} \to 0$ and $\|P_{t}\| \to \infty$, for $t \to \infty$. Thus, we consider this method practically useless.

One may consider it a serious downside of Taylor-approximation based methods that the gradient only approximates the shape of $f$ and thereby its mapping of the error on an `infinitesimally small neighborhood' of $m^-_{t+h}$. Hence, it might ignore the exact value of $y$ completely, if the mean prediction is far off. However, for a highly regular $f$  (e.g.~Lipschitz-continuous in the space variable) this gradient approximation is very good.

Moreover, the approximation by a first-order Taylor series expansion needs an approximation of the gradient, which explicit ODE solvers usually do not receive as an input. However, in many numerical algorithms (e.g.~optimization) the gradient is provided anyway. Therefore the gradient might already be known in real-world applications.
While we find this method promising when the gradient is known or can be efficiently computed, we exclude it from our experiments because the necessity of a gradient estimate breaks the usual framework of ODE solvers.

In contrast, Bayesian quadrature avoids the risk of a skewed distortion of the samples for Monte Carlo integration by actively spreading a grid of deterministic sigma-points. It does not need the gradient of $f$ and still can encode prior knowledge over $f$ by the choice of the covariance function if more is known \citep{Briol2015probint}. The potential of using Bayesian quadrature as a part of a filter was further explored by \citet{PS15}, however in the less structured setting of nonlinear filtering where additional inaccuracy from the linear approximation in the prediction step arises. Moreover, \citet{SHSV2015} recently pointed out that BQ can be seen as sigma-point methods and gave covariance functions and evaluation points which reproduce numerical integration methods known for their favorable behavior (for example Gauss--Hermite quadrature, which is used for a Gaussian weight function).

Due to these advantages, we propose a new class of BQ-based probabilistic ODE filters named BQ Filtering.

\subsection{BAYESIAN QUADRATURE FILTERING (BQF)}\label{BQ}

The crucial source of error for filtering-based ODE solvers is the calculation of the gradient measurement $y$ and its variance $R$ (c.f.~Section \ref{Framework}). We propose the novel approach to use BQ to account for the uncertainty of the input and thereby estimate $y$ and $R$. This gives rise a novel class of filtering-based solvers named \emph{BQ Filter} (BQF). As a filtering-based method, one BQF-step consists of the KF prediction step \eqref{KFP1}--\eqref{KFP2}, the calculation of $y$ and $R$ by BQ and the KF update step \eqref{KFU1}--\eqref{KFU5}. 

\begin{figure}[ht]
  \setlength{\figwidth}{.9\columnwidth}
  \setlength{\figheight}{5.75cm}
  \centering\scriptsize
  \input{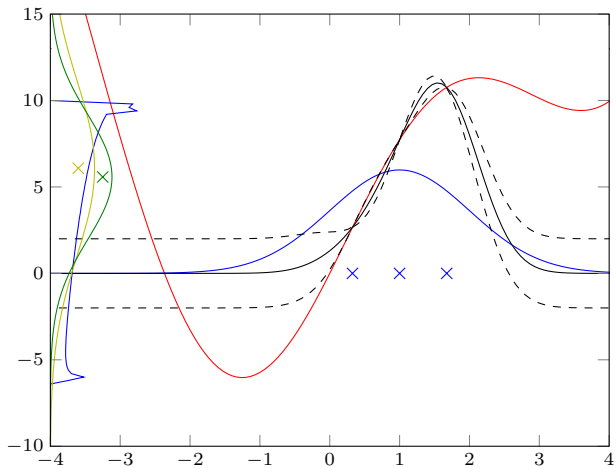}
  \caption{Prediction of function $f(x) = 8 \sin (x) + x^2$ (red) under uncertain input $x\sim \mathcal N(x;1,1)$ (density in blue). $\mathcal{GP}(f)$ (black) derived from Gaussian grid evaluation points with $N=3$ (blue crosses) as mean $\pm$ $2$ standard deviation. True distribution of prediction in blue. Gaussian fit to true distribution in yellow and predicted distribution by BQ in green with crosses at means.}
  \label{fig:pred}
\end{figure}

The KF prediction step outputs a Gaussian belief $\mathcal N(m_{t+h}^-,P_{t+h}^-)$ over the exact solution $u(t+h)$. This input value is propagated through $f$ yielding a distribution over the gradient at time $t+h$. In other words, our belief over $\nabla f(t+h,u(t+h))$ is equal to the distribution of $Y:=f(t,X)$, with uncertain input $X \sim \mathcal N(m_{t+h}^-,P_{t+h}^-)$. For general $f$ the distribution of $Y$ will be neither Gaussian nor unimodal (as e.g.~in Figure \ref{fig:pred}). But it is possible to compute the moments of this distribution under Gaussian assumptions on the input and the uncertainty over $f$ (see for example \citet{Deisenroth2009thesis}). The equivalent formulation of prediction under uncertainty clarifies as numerical integration clarifies the connection to sigma-point methods, i.e.~quadrature rules \citep{SHSV2015}. Quadrature is as extensively studied and well-understood as the solution of ODEs. A basic overview can be found in \citet{PTVF07}.
Marginalizing over $X$ yields an integral with Gaussian weight function
\begin{align}
  \E[Y] = \int f(t+h,x) \underbrace{\mathcal N(x;m_{t+h}^-,P^-_{t+h})}_{ =: w(x)} \ \rd x,
\end{align}
which is classically solved by quadrature, i.e.~evaluating $f$ at a number of evaluation points $(x_i)_{i=1,\dots,N}$ and calculating a weighted sum of these evaluations. BQ can be interpreted as a probabilistic extension of these quadrature rules in the sense that their posterior mean estimate of the integral coincides with classic quadrature rules, while adding a posterior variance estimate at low cost \citep{SHSV2015}. 

By choosing a kernel $k$ over the input space of $f$ and evaluation points $(x_i)_{i=1,\dots,N}$, the function $f$ is approximated by a GP regression \citep{RasmussenWilliams} with respect to the function evaluations $(f(x_i))_{i=1,\dots,N}$, yielding a GP posterior over $f$ with mean $m_f$ and covariance $k_f$ denoted by $\mathcal{GP}(f)$. The integral is then approximated by integrating the GP approximation, yielding the predictive distribution for $\mathcal I [f]$: 
\begin{align}\label{GPint}
  \mathcal I[f] \sim \int \mathcal{GP}(f)(x)\cdot \mathcal N(x;m_{t+h}^-,P^-_{t+h}) \ \rd x.
\end{align}

The uncertainty arising from the probability measure over the input is now split up in two parts: the uncertainty over the input value $x\sim \mathcal N(0,I)$ and the uncertainty over the precise value at this uncertain input, which can only be approximately inferred by its covariance with the evaluation points $(x_i)_{i=1,\dots,N}$, i.e.~by $\mathcal{GP}(f)$. These two kinds of uncertainty are depicted in Figure \ref{fig:pred}.
From the predictive distribution in \eqref{GPint}, we can now compute a posterior mean and variance of $\mathcal I[f]$ which results in a weighted sum for the mean
\begin{align}
  y:=\E\left[\mathcal I[f]\right] = \alpha^T K^{-1} \begin{pmatrix} f(x_1) , \hdots , f(x_n) \end{pmatrix}^T
\end{align}
with 
\begin{align}\label{alpha}
\alpha (i) = \int k(x,x_i) \mathcal N(x;0,I) \ \rd x
\end{align}
and variance
\begin{align}
  R&:=\Var \left[\mathcal I (f)\right] \notag \\ &= \int \int k(x,x^{\prime}) w(x) w(x^{\prime}) \ \rd x \rd x^{\prime} - \alpha^T K^{-1} \alpha,
\end{align}
where $K\in\R^{N\times N}$ denotes the kernel matrix, i.e.~$K_{i,j}=k(x_i,x_j)$.

The measurement generation in BQF is hence completely defined by the two free choices of BQ: the kernel $k$ and the evaluation points $(x_i)_{i=1,\hdots,n}$. By these choices, BQ and thereby the measurement generation in BQF is completely defined.
For the squared exponential kernel \citep{RasmussenWilliams}
\begin{align}\label{SEkernel}
  k(x,x^{\prime}) = \theta^2 \exp\left( -\frac 1{2\lambda^2} \|x-x^{\prime}\|^2  \right),
\end{align}
with lengthscale $\lambda>0$ and output variance $\theta^2>0$, it turns out that $y$ and $R$ can be computed in closed form and that many classic quadrature methods which are known for their favorable properties can be computed in closed form \citep{SHSV2015}, significantly speeding up computations. For the scalar case $nD=1$, we obtain for \eqref{alpha} by straightforward computations:
\begin{align}
  \alpha(i) = \frac{\lambda\theta^2}{\sqrt{\lambda^2+\sigma^2}} \exp\left(-\frac{(x_i-\mu)^2}{2(\lambda^2+\sigma^2)}\right), 
\end{align}
and
\begin{align}
\int \int k(x,x^{\prime}) w(x) w(x^{\prime}) \ \rd x \rd x^{\prime} = \frac{\theta^2}{\sqrt{1+2\sigma^2/\lambda^2}}
\end{align}

Hence, our BQ estimate for $y$ is given by the sigma-point rule
\begin{equation}
  y \approx \sum_{i=1}^N W_i f(t+h,x_i) 
\end{equation}
with easily computable weights
\begin{equation}
  W_i = [\alpha^T K^{-1}]_i.
\end{equation}
Also the variance $R$ takes a convenient shape
\begin{equation}
  R = \frac{\theta^2}{\sqrt{1+2\sigma^2/\lambda^2}} - \alpha^T K^{-1} \alpha.
\end{equation}
For $nD>1$, we get slightly more complicated formulas which are given in \citet{Deisenroth2009thesis}.

The other free choice in BQ, the evaluation points $(x_i)_{i=1,\dots,n}$, can also be chosen freely in every step of BQF. Usually, the nodes of BQ chosen are chosen so that the variance of the integral estimate is minimized (cf.~\citet{Briol2015probint}). For this algorithm, the uncertainty has to be measured, not minimized though. Hence, we propose just to take a uniform grid scaled by $\mathcal N(m_{t+h}^-,P_{t+h}^-)$ to measure the uncertainty in a comprehensive way. 

Another promising choice is given by the roots of the physicists' version of the Hermite polynomials, since they yield Gauss--Hermite quadrature (GHQ), the standard numerical quadrature against Gaussian measures, as a posterior mean for a suitable covariance function \citep{SHSV2015}. For GHQ, efficient algorithms to compute the roots and the weights are readily available \citep{PTVF07}.

\subsection{COMPUTATIONAL COST}\label{cost}

All of the presented algorithms buy their probabilistic extension to classic ODE solvers by adding computational cost, sometimes more sometimes less.  In most cases, evaluation of the vector field $f$ forms the computational bottleneck, so we will focus on it here. Of course, the internal computations of the solver adds cost as well. Since all the models discussed here have linear inference cost, though, this additional overhead is manageable. \\
The ML-algorithm by \citet{schober2014nips} is the fastest algorithm. By simply recasting a Runge--Kutta step as Gaussian filtering, rough probabilistic uncertainty is achieved with negligible computational overhead. \\
For the sampling method, the calculation of one individual sample of $Q_h$ amounts to running the entire underlying ODE solver once, hence the overall cost is $S$ times the original cost. \\ 
In contrast, the BQ-algorithm only has to run through $[0,T]$ once, but has to invert a $ND\times ND$ covariance matrix to perform Bayesian quadrature with $N$ evaluation points. Usually, $N$ will be small, since BQ performs well for a relatively small number of function evaluations (as e.g.~illustrated by the experiments below). However, if the output dimension $D$ is very large, Bayesian quadrature---like all quadrature methods---is not practical. BQ thus tends to be faster for small $D$, while MC tends to be faster for large $D$. 

When considering these computational overheads, there is a nuanced point to be made about the value-to-cost trade-off of constructing a posterior uncertainty measure. If a classic numerical solver of order $p$ is allotted a budget of $M$ times its original one, it can use it to reduce its step-size by a factor of $M$, and thus reduce its approximation error by an order $M^{p}$. It may thus seem pointless to invest even such a linear cost increase into constructing an uncertainty measure around the classic estimate. But, in some practical settings, it may be more helpful to have a notion of uncertainty on a slightly less precise estimate than to produce a more precise estimate without a notion of error. In addition, classic solvers are by nature sequential algorithms, while the probabilistic extensions (both the sampling-based and Gaussian-filtering based ones) can be easily parallelized. Where parallel hardware is available, the effective time cost of probabilistic functionality may thus be quite limited (although we do not investigate this possibility in our present experiments).

With regards to memory requirements, the MC-method needs significantly more storage, since it requires saving all sample paths, in order to statistically approximate the entire non-parametric measure $Q_h$ on $C^1([0,T],\R)$. The BQ-algorithm only has to save the posterior GP, i.e.~a mean and a covariance function, which is arguably the minimal amount to provide a notion of uncertainty. If MC reduces the approximation of $Q_h$ to its mean and variance, it only requires this minimal storage as well.

\section{EXPERIMENTS}\label{experiments}

\begin{figure*}[ht]\label{GroundTruth}
  \setlength{\figwidth}{1.8\columnwidth}
  \setlength{\figheight}{7.6cm}
  \centering\scriptsize
  \input{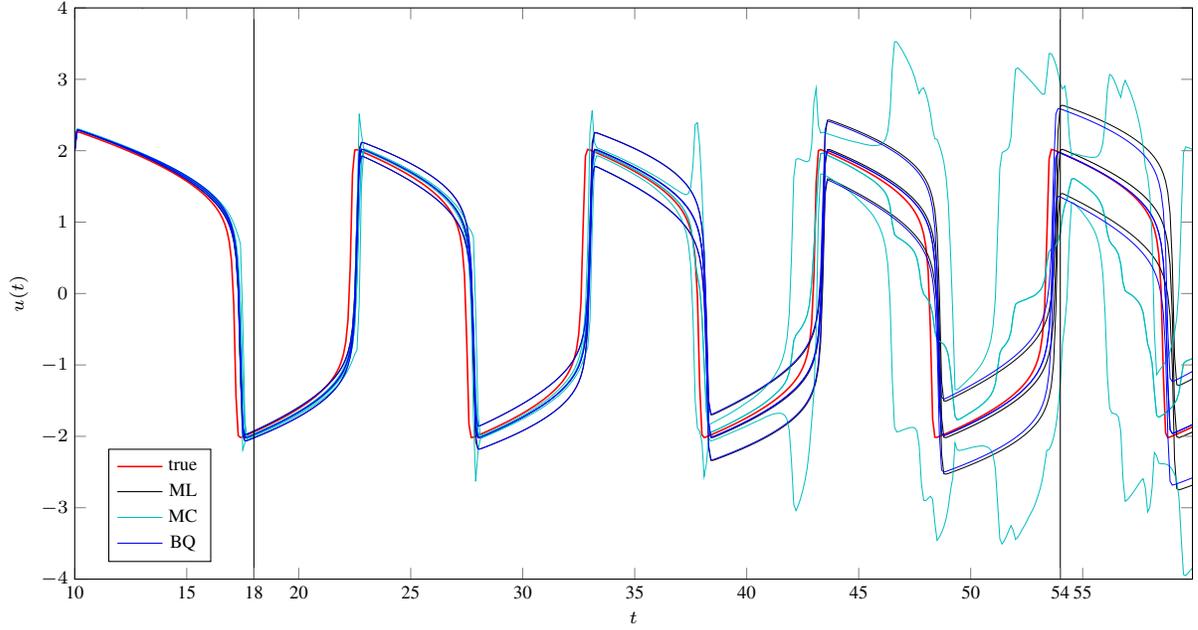}
  \caption{Solution estimates constructed on the Van der Pol oscillator \eqref{vdP}. True solution in red. Mean estimates of ML, MC and BQ in black, green, blue, respectively. Uncertainty measures (drawn at two times standard deviation) as thin lines of the same color.}
\end{figure*}

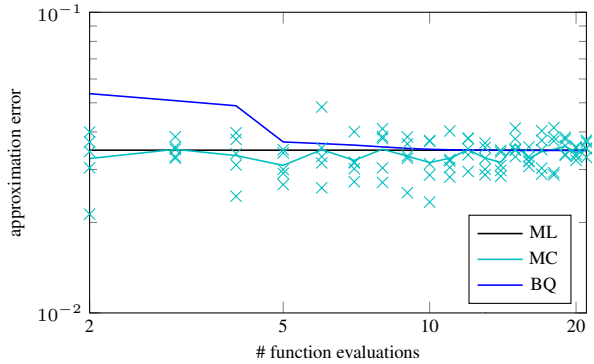
\begin{figure}[ht]\label{ErrorPlot1}
  \setlength{\figwidth}{.8\columnwidth}
  \setlength{\figheight}{4cm}
  \centering\scriptsize
  \input{exp909_Filter_fig0.tex}
  \caption{Plot of errors of the mean estimates at $t=18$ of the methods MC (green) and BQ (blue) as a function of the allowed function evaluations. Maximum likelihood error in black. Single runs of the probabilistic MC solver as green crosses. Average over all runs as green line.}
  \label{fig:exp909_0}
\end{figure}

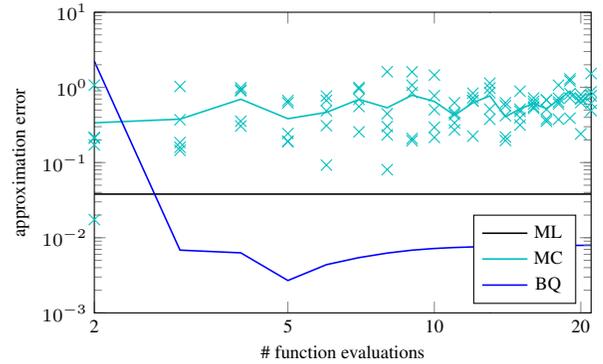
\begin{figure}[ht]\label{ErrorPlot2}
  \setlength{\figwidth}{.8\columnwidth}
  \setlength{\figheight}{4cm}
  \centering\scriptsize
  \input{exp909_Filter_fig1.tex}
  \caption{Plot of errors of the mean estimates at $t=54$ of the methods MC (green) and BQ (blue) as a function of the allowed function evaluations. Maximum likelihood error in black. Single runs of the probabilistic MC solver as green crosses. Average over all runs as green line.}
  \label{fig:exp909_1}
\end{figure}

This section explores applications of the probabilistic ODE solvers discussed in Section \ref{Background}. The sampling-based algorithm by \citep{conrad_old_arxiv} will be abbreviated as MC, the maximum-likelihood Gaussian filter (\citep{schober2014nips}) as ML and our novel BQ-based filter (BQF) as BQ. In particular, we assess how the performance of the purely deterministic class of Gaussian filtering based solvers compares to the inherently random class of sampling-based solvers.

We experiment on the Van der Pol oscillator \citep{hairer87:_solvin_ordin_differ_equat_i}, a non-conservative oscillator with non-linear damping, which is a standard example for a non-stiff dynamical system. It is governed by the equation
\begin{equation} \label{vdP}
  \frac{\partial^2 u}{\partial t^2}  = \mu (1-u^2)\frac{\partial u}{\partial t} - u,
\end{equation}
where the parameter $\mu\in\R$ indicates the non-linearity and the strength of the damping. We set $\mu = 5$ on a time axis $[10,60]$, with initial values $(u(10),\dot{u}(10)) = (2,10)$. 

All compared methods use a model of order $q=3$, and a step size $h=0.01$. This induces a state-space model given by a twice-integrated Wiener process prior (cf.~\eqref{SDE}) which yields a version of ML close to second-order Runge--Kutta \citep{schober2014nips}. The same solver is used as the underlying numerical solver $\Psi_h$ in MC. For the noise parameter, which scales the deviation of the evaluation point of $f$ from the numerical extrapolation (i.e.~the variance of the driving Wiener process for ML and BQ, and the variance of $\xi_k$ for MC), we choose $\sigma^2 = 0.1$. The drift matrix $F$ of the underlying integrated Wiener process is set to the default values $f_i = i$ for $i=1,\dots,q-1$. The covariance function used in BQ is the widely popular squared exponential \eqref{SEkernel}, with lengthscale $\lambda = 1$ and output variance $\theta^2 = 1$. (Since all methods use the same model, this tuning does not favor one algorithm over the other. In practice all these parameters should of course be set by statistical estimation.).

For a fair comparison in all experiments, we allow MC and BQ to make the same amount of function evaluations per time step. If MC draws $N$ samples, BQ uses $N$ evaluation points.
The first experiment presents the solutions of the presented algorithms on the van der Pol oscillator \eqref{vdP} on the whole time axis in one plot, when we allow BQ and MC to make five function evaluations.
Then, we examine more closely how the error of each methods changes as a function of the number of evaluations of $f$ in Figure \ref{fig:exp909_0} and Figure \ref{fig:exp909_1}.

\subsection{SOLUTION MEASURES ON VAN DER POL OSCILLATOR}

Figure \ref{GroundTruth} shows the solution estimates constructed by the three solvers across the time domain. In all cases, the mean estimates roughly follow the exact solution (which e.g.~Gaussian filtering with Monte Carlo integration by sampling \eqref{MCIS1}--\eqref{MCIS3} does not achieve). A fundamental difference between the filtering-based methods (ML and BQ) and the sampling-based MC algorithm is evident in both the mean and the uncertainty estimate.

While the filtering-based methods output a trajectory quite similar to the exact solution with a small time lag, the MC algorithm produces a trajectory of a more varying shape. Characteristic points of the MC mean estimate (such as local extrema) are placed further away from the exact value than for filtering-based methods.

The uncertainty estimation of MC appears more flexible as well. ML and BQ produce an uncertainty estimate which runs parallel to the mean estimate and appears to be strictly increasing. It appears to increase slightly in every step, resulting in an uncertainty estimate, which only changes very slowly. The solver accordingly appears overconfident in the spikes and underconfident in the valleys of the trajectory. The uncertainty of MC varies more, scaling up at the steep parts of the oscillator and decreasing again at the flat parts, which is a desirable feature.

Among the class of filtering-based solvers, the more refined BQ method outputs a better mean estimate with more confidence than ML.

\subsection{QUALITY OF ESTIMATE AS A FUNCTION OF ALLOWED EVALUATIONS}

Figure \ref{fig:exp909_0} and Figure \ref{fig:exp909_1} depict the value of the error of the mean approximation as a function of the allowed function evaluations $N$ (i.e.~$N$ evaluation points for BQ and $N$ samples for MC) at time points $t_1=18$ and $t_2 = 54$. Since the desired solution measure $Q_h$ for MC can only be statistically approximated by the $N$ samples, the mean estimate of MC is random. For comparison, the average of five MC-runs is computed.

At the early time point $t_1=18$, all trajectories are still close together and the methods perform roughly the same, as we allow more evaluations. There is a slight improvement for BQ with more evaluations, but the error remains above the one of ML error.

At the later time $t_2=54$, BQ improves drastically when at least five evaluations are allowed, dropping much below the ML error. 

The average error by MC appears to be not affected by the number of samples. The ML error is constant, because it always evaluates only once.

\section{DISCUSSION}

The conducted experiments provide an interesting basis to discuss the differences between filtering-based methods (ML and BQ) and the sampling-based MC algorithm. We make the following observations:

\begin{enumerate}[(i)]

\item \emph{Additional samples do not improve the random mean estimate of MC in expectation:}\\
Since the samples of MC are independent and identically distributed, the expectation of the random mean estimate of MC is the same, regardless of the amounts of samples. This property is reflected in Figure \ref{fig:exp909_0} and Figure \ref{fig:exp909_1}, by the constant green line (up to random fluctuation). Additional samples are therefore only useful to improve the uncertainty calibration.

\item \emph{The uncertainty calibration of MC appears more adaptive than of ML and BQ:}\\ Figure \ref{GroundTruth} suggests that MC captures the uncertainty more flexibly: It appropriately scales up in the steep parts of the oscillator, while expressing high confidence in the flat parts of the oscillator. The exact trajectory is inside the interval between mean $\pm$ $2$ standard deviations, which is not the case for BQ and ML. Moreover, MC produces a more versatile measure. The filtering-based methods appear to produce a strictly increasing uncertainty measure by adding to the posterior uncertainty in every step. MC avoids this problem by sampling multiple time over the whole time interval. We deem the resulting flexibility a highly desirable feature. BQ also outputs a meaningful uncertainty measure and we expect that adding Bayesian smoothing \citep{Sarkka2013} would enable filtering-based methods to produce more adaptive measures as well.

\item \emph{The expected error of MC-samples (and their mean) is higher than the error of ML:}\\ In the experiments, MC produced a higher error for the mean estimate, compared to both ML and BQ. We expect that this happens on all dynamical systems \emph{by construction}: Given $U_k$, the next value $U_{k+1}$ of a MC-sample is calculated by adding Gaussian noise $\xi_k$ to the ML-extrapolation starting in $U_k$ (cf.~equation \eqref{Conrad_algo}). Due to the symmetry and full support of Gaussian distributions, the perturbed solution has a higher error than the unperturbed prediction, which coincides with the ML solution. Hence, every MC-sample accumulates with every step a positive expected error increment compared to the ML estimate. By the linearity of the average, the mean over all samples inherits the same higher error than the ML mean (and thereby also than the error of the more refined BQ mean).
\end{enumerate}

Summing up, we argue that---at their current state---filtering-based methods appear to produce a `better' mean estimate, while sampling-based methods produce in some sense a `better' uncertainty estimate. Many applications might put emphasis on a good mean estimate, while needing a still well-calibrated uncertainty quantification. Our method BQF provides a way of combining a precise mean estimate with a meaningful uncertainty calibration. Sampling-based methods might not be able to provide this due to their less accurate mean estimate. For future work (which is beyond the scope of this paper), it could be possible to combine the advantages of both approaches in a unified method.

\section{CONCLUSION}

We have presented theory and methods for the probabilistic solution of ODEs which provide uncertainty measures over the solution of the ODE, contrasting the classes of (deterministic) filtering-based and (random) sampling-based solvers. We have provided a theoretical framework for Gaussian filtering as state space inference in linear Gaussian SDEs, highlighting the prediction of the gradient as the primary source of uncertainty. Of all investigated approximations of the gradient, Bayesian Quadrature (BQ) produces the best results, by actively learning the shape of the vector field $f$ through deterministic evaluations. Hence, we propose a novel filtering-based method named \emph{Bayesian Quadrature Filtering} (BQF), which employs BQ for the gradient measurement.  

For the same amount of allowed gradient evaluations, the mean estimate of BQF appears to outperform the mean estimate of state-of-the-art sampling-based solvers on the Van der Pol oscillator, while outputting a better calibrated uncertainty than other filtering-based methods.


\bibliography{../../../../../../bibfile.bib}
\bibliographystyle{icml2016}

\end{document}

%% file: mypreamble.tex
\usepackage{amsmath,amssymb}
\usepackage{pgfplots}
\usepackage[]{amsmath,amssymb,amsthm}
\usepackage{paralist}
\usepackage[applemac]{inputenc}
\usepackage{graphicx}
\usepackage{setspace}
\usepackage{enumerate}
\usepackage{natbib}
\usepackage{bbm}
\usepackage{pgfplots}
\DeclareMathOperator{\Var}{Var}
\pgfplotsset{compat=newest}
\newlength{\figwidth}
\newlength{\figheight}
\usetikzlibrary{external}

\newcommand{\R}{\mathbb{R}}

\newcommand{\N}{\mathbb{N}}

\newcommand{\E}{\mathbb{E}}

\newcommand{\del}{\delta}
\newcommand{\vertiii}[1]{{\left\vert\kern-0.25ex\left\vert\kern-0.25ex\left\vert #1 
   \right\vert\kern-0.25ex\right\vert\kern-0.25ex\right\vert}} 
\usepackage{mathtools}

%% file: exp909_Filter_fig0.tex
%
%
%
%
\begin{tikzpicture}

\definecolor{color0}{rgb}{0,0.75,0.75}

\begin{axis}[
xmin=2, xmax=21,
ymin=0.01, ymax=0.1,
xmode=log,
ymode=log,
axis on top,
width=\figwidth,
height=\figheight,
ylabel near ticks,
xlabel near ticks,
xlabel={\# function evaluations},
ylabel={approximation error},
scale only axis,
xtick={1, 2, 5, 10, 20},
xticklabels={1,2,5,10,20},
legend entries={{ML},{MC},{BQ}},
legend style={at={(0.97,0.03)}, anchor=south east}
]
\addplot [semithick, black]
coordinates {
(2,0.0347819934753364)
(3,0.0347819934753364)
(4,0.0347819934753364)
(5,0.0347819934753364)
(6,0.0347819934753364)
(7,0.0347819934753364)
(8,0.0347819934753364)
(9.00000000000001,0.0347819934753364)
(10,0.0347819934753364)
(11,0.0347819934753364)
(12,0.0347819934753364)
(13,0.0347819934753364)
(14,0.0347819934753364)
(15,0.0347819934753364)
(16,0.0347819934753364)
(17,0.0347819934753364)
(18,0.0347819934753364)
(19,0.0347819934753364)
(20,0.0347819934753364)
(21,0.0347819934753364)

};
\addplot [semithick, color0]
coordinates {
(2,0.032619134382048)
(3,0.0350074953610477)
(4,0.0333719379419939)
(5,0.030963085261862)
(6,0.0349082815907346)
(7,0.0321704812450629)
(8,0.0350244979885586)
(9.00000000000001,0.0332179519309591)
(10,0.0316339218892346)
(11,0.0326939933285086)
(12,0.0348877194413803)
(13,0.0327593336836235)
(14,0.0316507964076719)
(15,0.0364081370056302)
(16,0.0326668510820485)
(17,0.0344030708110202)
(18,0.0350191513081666)
(19,0.0358649628596575)
(20,0.0340274284247553)
(21,0.0355863185442107)

};
\addplot [semithick, blue]
coordinates {
(2,0.0536162903987345)
(3,0.0508254977848808)
(4,0.0488710480182026)
(5,0.0370436933557461)
(6,0.0365784198440753)
(7,0.036159984017339)
(8,0.0356670829577808)
(9.00000000000001,0.0352591304601788)
(10,0.0350189532654588)
(11,0.0348933163200873)
(12,0.0348259089727638)
(13,0.0347867571259914)
(14,0.0347618684937792)
(15,0.0347446499988637)
(16,0.034731830318439)
(17,0.0347216833808099)
(18,0.0347132433755375)
(19,0.0347059425483662)
(20,0.0346994351148115)
(21,0.0346935067765828)

};
\addplot [color0, mark=x, mark size=3, only marks]
coordinates {
(2,0.0213050079973465)
(3,0.0326978996957525)
(4,0.0338512468887071)
(5,0.0267253907720053)
(6,0.0483888096335659)
(7,0.0301803046450224)
(8,0.040984971101238)
(9.00000000000001,0.0384997130003952)
(10,0.0313570666443257)
(11,0.0282866103742443)
(12,0.035118008711231)
(13,0.0355068084098882)
(14,0.0342583425989018)
(15,0.0412125398897107)
(16,0.0330963945387668)
(17,0.0303278984200626)
(18,0.0292792711225738)
(19,0.0359807096080895)
(20,0.0331653726254597)
(21,0.0370053680342799)

};
\addplot [color0, mark=x, mark size=3, only marks]
coordinates {
(2,0.0370390095821484)
(3,0.0385623187918098)
(4,0.0244402536021624)
(5,0.0348862886237575)
(6,0.0315488544013414)
(7,0.0315361785560935)
(8,0.0271668401199194)
(9.00000000000001,0.0332726204879061)
(10,0.0372404841128027)
(11,0.0324223391632787)
(12,0.0380410227352332)
(13,0.0288148791473561)
(14,0.0309649285635181)
(15,0.0379625795023453)
(16,0.0306979607611522)
(17,0.0403298789111943)
(18,0.0383631136304749)
(19,0.0336747184572312)
(20,0.0321682537136385)
(21,0.0352702739322774)

};
\addplot [color0, mark=x, mark size=3, only marks]
coordinates {
(2,0.0343987340335217)
(3,0.0350465456740923)
(4,0.0397498791611031)
(5,0.02913066565822)
(6,0.0331249926705468)
(7,0.0273157761542762)
(8,0.0303568967861489)
(9.00000000000001,0.0251119677285423)
(10,0.0374453473280136)
(11,0.0304718925202281)
(12,0.0335776861596044)
(13,0.0367511652646724)
(14,0.0347632156360935)
(15,0.0347685583658337)
(16,0.0318493523730825)
(17,0.0289401601025123)
(18,0.0373313830026443)
(19,0.0334360250320607)
(20,0.0357278129101184)
(21,0.035177486767408)

};
\addplot [color0, mark=x, mark size=3, only marks]
coordinates {
(2,0.0303944366052748)
(3,0.0331160583603345)
(4,0.0377732310873005)
(5,0.0340707531073234)
(6,0.026067442918978)
(7,0.0317514381697279)
(8,0.0380341165683702)
(9.00000000000001,0.0326790997570019)
(10,0.0287524505839656)
(11,0.0320090644600159)
(12,0.0295425272104495)
(13,0.032679041646799)
(14,0.0297892991262121)
(15,0.0350077280445331)
(16,0.0320234609032317)
(17,0.0375123338908931)
(18,0.0413438508262713)
(19,0.0378172686448628)
(20,0.0337364314295203)
(21,0.0375503200123472)

};
\addplot [color0, mark=x, mark size=3, only marks]
coordinates {
(2,0.0399584836919487)
(3,0.0356146542832494)
(4,0.0310450789706966)
(5,0.030002328148004)
(6,0.0354113083292409)
(7,0.0400687087001947)
(8,0.0385796653671162)
(9.00000000000001,0.0365263586809499)
(10,0.0233742607770653)
(11,0.0402800601247759)
(12,0.0381593523903831)
(13,0.0300447739494016)
(14,0.0284781961136342)
(15,0.0330892792257282)
(16,0.0356670868340094)
(17,0.0349050827304387)
(18,0.0287781379588685)
(19,0.0384160925560431)
(20,0.0353392714450393)
(21,0.0329281439747411)

};

\path [draw=black, fill opacity=0] (axis cs:13,0.1)--(axis cs:13,0.1);

\path [draw=black, fill opacity=0] (axis cs:21,13)--(axis cs:21,13);

\path [draw=black, fill opacity=0] (axis cs:13,0.01)--(axis cs:13,0.01);

\path [draw=black, fill opacity=0] (axis cs:2,13)--(axis cs:2,13);

\end{axis}

\end{tikzpicture}

%% file: exp909_Filter_fig1.tex
%
%
%
%
\begin{tikzpicture}

\definecolor{color0}{rgb}{0,0.75,0.75}

\begin{axis}[
xmin=2, xmax=21,
ymin=0.001, ymax=10,
xmode=log,
ymode=log,
axis on top,
width=\figwidth,
height=\figheight,
ylabel near ticks,
xlabel near ticks,
xlabel={\# function evaluations},
ylabel={approximation error},
scale only axis,
xtick={1, 2, 5, 10, 20},
xticklabels={1,2,5,10,20},
legend entries={{ML},{MC},{BQ}},
legend style={at={(0.97,0.03)}, anchor=south east}
]
\addplot [semithick, black]
coordinates {
(2,0.038065355104439)
(3,0.038065355104439)
(4,0.038065355104439)
(5,0.038065355104439)
(6,0.038065355104439)
(7,0.038065355104439)
(8,0.038065355104439)
(9.00000000000001,0.038065355104439)
(10,0.038065355104439)
(11,0.038065355104439)
(12,0.038065355104439)
(13,0.038065355104439)
(14,0.038065355104439)
(15,0.038065355104439)
(16,0.038065355104439)
(17,0.038065355104439)
(18,0.038065355104439)
(19,0.038065355104439)
(20,0.038065355104439)
(21,0.038065355104439)

};
\addplot [semithick, color0]
coordinates {
(2,0.338101423524576)
(3,0.378082733965294)
(4,0.695642859906598)
(5,0.383117763602541)
(6,0.465483608689959)
(7,0.692479490360809)
(8,0.535853676702201)
(9.00000000000001,0.786954650035278)
(10,0.649810719139642)
(11,0.425530736823475)
(12,0.623404727344167)
(13,0.771600035043177)
(14,0.406502230057445)
(15,0.518660974942948)
(16,0.611527228009649)
(17,0.520087788659293)
(18,0.681533291866716)
(19,0.890691376678115)
(20,0.617661027945417)
(21,0.856555769561829)

};
\addplot [semithick, blue]
coordinates {
(2,2.22854856713055)
(3,0.00683512621670701)
(4,0.00629361439387121)
(5,0.00269648756263208)
(6,0.00437909340349551)
(7,0.00543774797281116)
(8,0.00624760283290193)
(9.00000000000001,0.00683320737281967)
(10,0.00719428276718603)
(11,0.00740982542043822)
(12,0.00754601571341462)
(13,0.00763885748406512)
(14,0.0077065462449839)
(15,0.00775860234796233)
(16,0.00780032836771882)
(17,0.00783487364886026)
(18,0.00786421996579456)
(19,0.00788967642300208)
(20,0.00791214168172138)
(21,0.00793225053120716)

};
\addplot [color0, mark=x, mark size=3, only marks]
coordinates {
(2,0.2112518156424)
(3,0.371265231483039)
(4,0.899240719454549)
(5,0.189371798272776)
(6,0.76631008801949)
(7,0.25637682636375)
(8,1.61725928531763)
(9.00000000000001,1.60938627320202)
(10,0.504621866083098)
(11,0.459569776542358)
(12,0.664948466229009)
(13,1.13783017032189)
(14,0.219098007473559)
(15,0.480508607932043)
(16,0.552852888261122)
(17,0.353888533918339)
(18,0.6416347586758)
(19,0.389393924162856)
(20,0.716108266443921)
(21,0.491458684192163)

};
\addplot [color0, mark=x, mark size=3, only marks]
coordinates {
(2,0.217262895708699)
(3,0.159633000348915)
(4,0.355007082390613)
(5,0.667681167901081)
(6,0.309738831845013)
(7,0.954372803585452)
(8,0.446259871385954)
(9.00000000000001,0.850411711157405)
(10,0.776472670107834)
(11,0.271296796459088)
(12,0.224311817842782)
(13,0.376975879922057)
(14,0.196683527237885)
(15,0.394173667461424)
(16,0.653036507218052)
(17,0.680275833721228)
(18,1.07014589944359)
(19,0.82553482914524)
(20,0.240953461957615)
(21,0.615125313786819)

};
\addplot [color0, mark=x, mark size=3, only marks]
coordinates {
(2,1.07270859545041)
(3,0.145502446582283)
(4,0.989622691561405)
(5,0.243359026018282)
(6,0.67114833692894)
(7,1.00041906137748)
(8,0.301744094648146)
(9.00000000000001,0.195288829902295)
(10,1.45707122191948)
(11,0.337160624111709)
(12,0.837727632657531)
(13,0.848794947627503)
(14,0.420789041956778)
(15,0.890071946931995)
(16,0.60382197552142)
(17,0.509383965717388)
(18,0.375831047481676)
(19,0.724277633846242)
(20,0.629094820285667)
(21,0.875978531145223)

};
\addplot [color0, mark=x, mark size=3, only marks]
coordinates {
(2,0.17185845585199)
(3,1.02927279582539)
(4,0.928289922349506)
(5,0.192874365040189)
(6,0.487660928087234)
(7,0.694297290733424)
(8,0.0804954598067338)
(9.00000000000001,0.209350125130936)
(10,0.294703329834488)
(11,0.624701083157881)
(12,0.662605018370608)
(13,0.513740678884718)
(14,0.624330378896481)
(15,0.315844309233602)
(16,0.526367935339026)
(17,0.679041135321522)
(18,0.71877139438785)
(19,1.29522452379368)
(20,0.696933772750579)
(21,1.52676441287424)

};
\addplot [color0, mark=x, mark size=3, only marks]
coordinates {
(2,0.0174253549693826)
(3,0.184740195586841)
(4,0.306053883776918)
(5,0.622302460780376)
(6,0.0925598585691189)
(7,0.556931469743936)
(8,0.233509672352543)
(9.00000000000001,1.07033631078373)
(10,0.216184507753315)
(11,0.43492540384634)
(12,0.727430701620904)
(13,0.980658498459714)
(14,0.571610194722524)
(15,0.512706343155673)
(16,0.721556833708627)
(17,0.377849474617988)
(18,0.601283359344661)
(19,1.21902597244256)
(20,0.805214818289299)
(21,0.773451905810697)

};

\path [draw=black, fill opacity=0] (axis cs:13,10)--(axis cs:13,10);

\path [draw=black, fill opacity=0] (axis cs:21,13)--(axis cs:21,13);

\path [draw=black, fill opacity=0] (axis cs:13,0.001)--(axis cs:13,0.001);

\path [draw=black, fill opacity=0] (axis cs:2,13)--(axis cs:2,13);

\end{axis}

\end{tikzpicture}